\documentclass[reqno]{amsart}
\usepackage{amsmath, amssymb,amsthm,bm,mathrsfs}
\usepackage[utf8]{inputenc}
\usepackage{enumerate}
\usepackage{hyperref}
\usepackage{graphicx}
\usepackage{pgf,tikz,pgfplots}
\pgfplotsset{compat=1.16}
\usetikzlibrary{arrows.meta,patterns,matrix,decorations.pathreplacing,calc,positioning}

\usepackage{todonotes}

\newtheorem{theorem}{Theorem}

\theoremstyle{definition}
\newtheorem{definition}{Definition}

\theoremstyle{remark}
\newtheorem{remark}{Remark}

\title[A conditional bound on sphere tangencies]{A conditional bound on sphere tangencies in all dimensions}
\author[C. Crowley]{Conrad Crowley}
\address[C. Crowley]{School of Mathematical Sciences, University College Cork, Western Gateway Building, Western Road, Cork, Ireland}
\email{conradcrowley@gmail.com}

\author[M. Vitturi]{Marco Vitturi}
\address[M. Vitturi]{School of Mathematical Sciences, University College Cork, Western Gateway Building, Western Road, Cork, Ireland}
\email{marco.vitturi@ucc.ie}

\begin{document}
\begin{abstract}
We use polynomial method techniques to bound the number of tangent pairs in a collection of $N$ spheres in $\mathbb{R}^n$ subject to a non-degeneracy condition, for any $n \geq 3$. The condition, inspired by work of Zahl for $n=3$, asserts that on any sphere of the collection one cannot have more than $B$ points of tangency concentrated on any low-degree subvariety of the sphere. For collections that satisfy this condition, we show that the number of tangent pairs is $O_{\epsilon}(B^{1/n - \epsilon} N^{2 - 1/n + \epsilon})$.
\end{abstract}

\maketitle
%\tableofcontents

% Create a new 1st level heading
\section{Introduction}
\begin{figure}[t]
\centering
\includegraphics[scale=0.2]{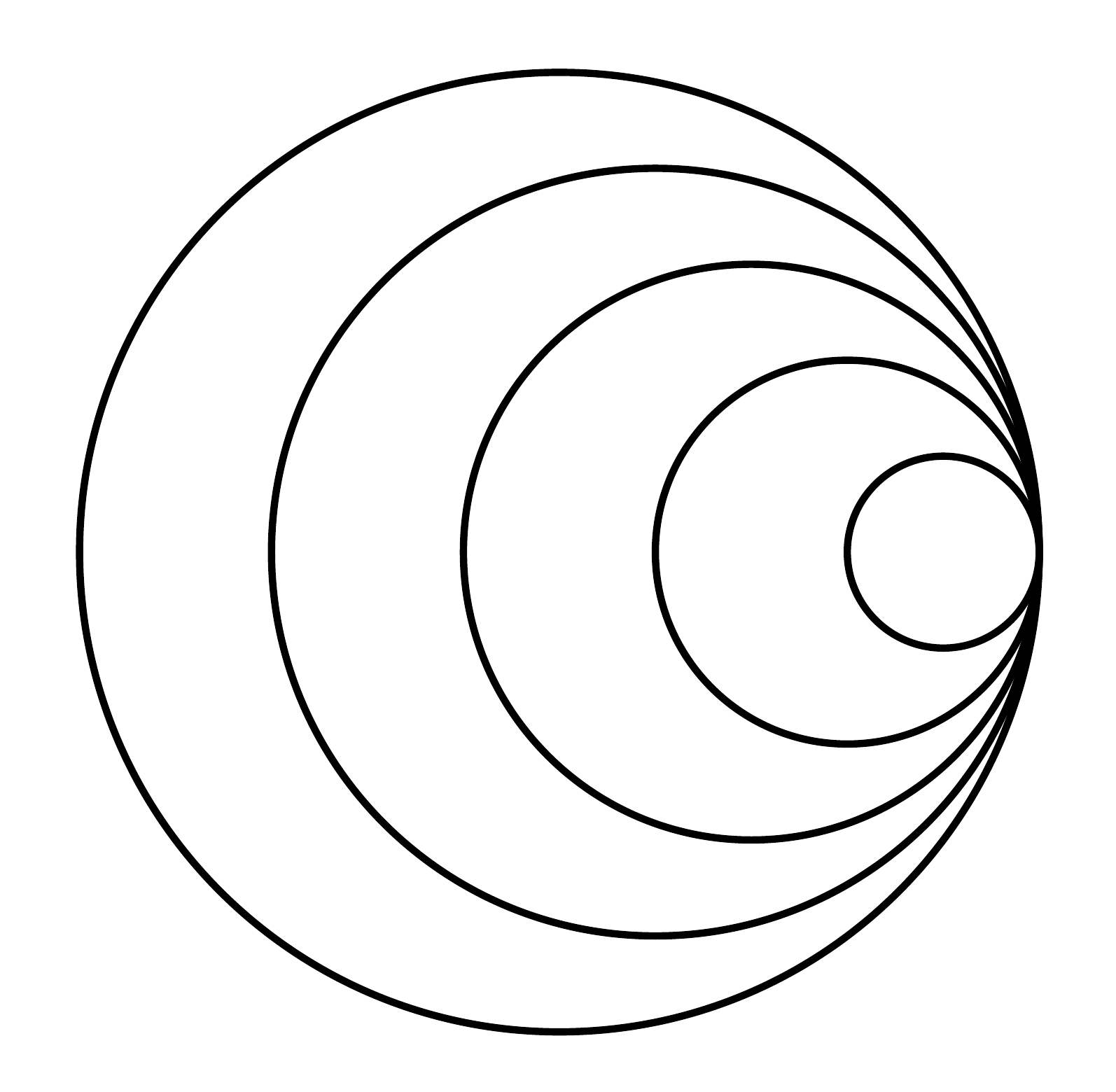}
\caption{\footnotesize A configuration in which all spheres are pairwise tangent.} \label{figure:hawaiian_ring}
\end{figure}

Let $n\geq 2$ and consider a finite collection $\mathscr{C}$ of spheres of arbitrary radii and centres in $\mathbb{R}^n$. Let $N$ be the cardinality of $\mathscr{C}$. An interesting problem in combinatorics is to count the maximum number of pairs of tangent spheres in the collection $\mathscr{C}$. It is always possible to arrange the $N$ spheres so that they are all pairwise tangent to each other at a single point (see Figure \ref{figure:hawaiian_ring}), so that in order to obtain a bound for the number of tangent pairs which is asymptotically better than $O(N^2)$ some non-degeneracy conditions must be assumed.\par
When $n=2$ this problem was considered first by Wolff in \cite{Wolff} and again by Ellenberg, Solymosi and Zahl in \cite{EllenbergSolymosiZahl} (the latter also for general bounded degree curves and fields other than $\mathbb{R}$). They imposed the non-degeneracy condition that no three circles in $\mathscr{C}$ can be mutually tangent to each other at a same point; under this condition, in \cite{EllenbergSolymosiZahl} it was shown that the number of tangent pairs is $O(N^{3/2})$.\footnote{Wolff obtained $O_{\epsilon}(N^{3/2+\epsilon})$ in \cite{Wolff} instead.} An alternative proof of the same bound can be found in \cite{SharirSmorodinskyValculescudeZeeuw}. It is not known whether the exponent $3/2$ is asymptotically sharp -- in upcoming work, Pohoata has shown that the number of tangent pairs is at most $o(N^{3/2})$ (personal communication). The best lowerbound available on the other hand is $\gtrsim N^{4/3}$ (obtained via inversion from a configuration of $N/2$ lines and $N/2$ points in $\mathbb{R}^2$ which are critical for the Szemer\'{e}di-Trotter inequality; see Section 8.1 of \cite{EllenbergSolymosiZahl}). The proof in \cite{EllenbergSolymosiZahl} is a combination of a lifting procedure that translates the tangency problem into a problem of incidence geometry for curves in $\mathbb{R}^3$ and of a polynomial method akin to that of Dvir from \cite{Dvir}. The proof in \cite{SharirSmorodinskyValculescudeZeeuw} is based on a different lifting process that translates the problem into one of point-surface incidences, for which one can apply the general incidence bounds of \cite{FoxPachShefferSukZahl}.\par
\begin{figure}[h]
\centering
\includegraphics[scale=0.3]{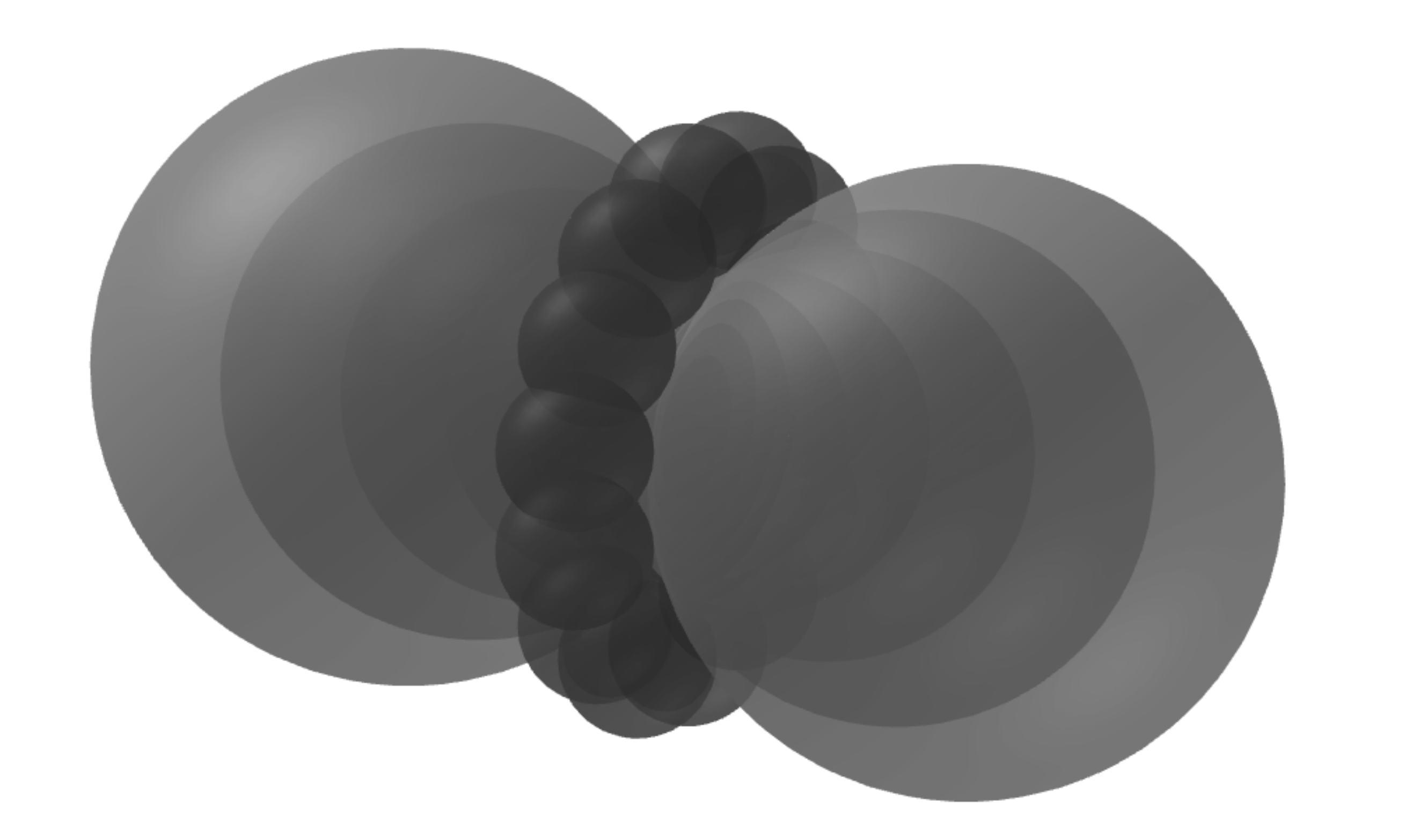}
\caption{\footnotesize A pair of complementary conic sections as per \cite{Zahl}, that is, a configuration consisting of two subcollections such that every sphere is tangent to every other sphere in the other subcollection.} \label{figure:pair_of_conics}
\end{figure}

When $n=3$ the problem was considered by Zahl in \cite{Zahl}, in which he addressed the more general situation of spheres with signed radii. Signed radii allow to distinguish internal from external tangencies: if $c_1, c_2$ are the centres of two spheres of signed radii $r_1, r_2$, the spheres are in contact if $\|c_1 - c_2\|^2 = (r_1 - r_2)^2$, and more specifically are internally tangent if $r_1, r_2$ have the same sign, and externally tangent otherwise. In this dimension there are additional configurations that need to be ruled out in order to obtain a non-trivial bound, besides those of Figure \ref{figure:hawaiian_ring}. Indeed, there exist collections of spheres consisting of two subcollections $\mathscr{C}', \mathscr{C}''$ such that any sphere in $\mathscr{C}'$ is tangent to every sphere in $\mathscr{C}''$ and viceversa; such collections can therefore contain up to $\sim N^2$ tangent pairs. Such a configuration, referred to in \cite{Zahl} as a pair of complementary conic sections,\footnote{If we identify the sphere of centre $(x,y,z)$ and radius $r$ with the coordinates $(x,y,z,r)$, then indeed the spheres in such configurations correspond to points on a pair of conic sections of $\mathbb{R}^4$.} is illustrated in Figure \ref{figure:pair_of_conics}. Zahl was able to show that in this dimension these are essentially the only problematic configurations, and it follows from his results that if the collection $\mathscr{C}$ is such that 
\begin{enumerate}[(i)]
\item no three spheres are mutually tangent to each other at a same point,
\item any pair of complementary conic sections contains at most $N^{1/2}$ spheres,
\end{enumerate}
then the collection contains at most $O(N^{3/2})$ pairs of tangent spheres.  The proof relies on an elegant fact known as the Lie sphere-line correspondence, which allows to translate the tangency problem into one of incidence geometry for complex lines in a copy of the Heisenberg group embedded in $\mathbb{C}^3$. The exponent $3/2$ is sharp, as can be seen by taking a collections of spheres with centres $(a,b,c)$ and radii $d$ where $a,b,c,d \in \{1, \ldots, M\}$ (see Example 1.17 of \cite{Zahl}).\par
In dimensions higher than $3$ it appears likely that additional configurations would need to be ruled out in order to obtain a non-trivial bound on the number of tangent pairs; the number of such configurations is also potentially dependent on the dimension. To sidestep this issue we will introduce a non-degeneracy condition that takes the same form in every dimension; this will allow us to obtain bounds that beat the trivial $O(N^2)$ one (conditionally on a parameter). To motivate the condition we are about to introduce, observe that on any given sphere in the configurations presented in Figure \ref{figure:pair_of_conics} the points of tangency all lie on the intersection of the sphere with a plane -- that is, there are $\sim N$ points of tangency concentrated on a low-degree sub-variety of the sphere. Inspired by this observation, we make the following definition.
\begin{definition}
Let $\mathscr{C}$ be a finite collection of spheres in $\mathbb{R}^n$, and assume that no three spheres are mutually tangent to each other at a same point. For $B,D>0$, we say that $\mathscr{C}$ is \emph{$(B,D)$-non-degenerate} if for every sphere $\sigma \in \mathscr{C}$ and every sub-variety $\Gamma \subset \sigma$ of degree at most $2D$ we have that the number of points of $\Gamma$ at which $\sigma$ is tangent to another sphere in $\mathscr{C}$ is at most $B$.
\end{definition}
Heuristically, the condition can be understood as a statement about the ``algebraic complexity'' of the set of tangencies: on any given sphere in a non-degenerate collection, the points of tangency cannot cluster on low-degree subvarieties; that is to say, to interpolate many points of tangency at once a high-degree polynomial is needed.\par
For collections of $N$ spheres that satisfy such a condition, we are able to show that the number of tangent pairs is asymptotically smaller than $O(N^2)$, provided $B$ is much smaller than $N$. In particular, our main result is the following.
\begin{theorem}\label{main_theorem}
Let $n\geq 3$. For any $\epsilon > 0$ there exists a degree $D = D(\epsilon,n)$ such that the following holds. Let $\mathscr{C}$ be a collection of $N$ spheres in $\mathbb{R}^n$ and let $B \leq N$. If $\mathscr{C}$ is $(B,D)$-non-degenerate, we have 
\[ |\{ (\sigma,\sigma') \in \mathscr{C} \times \mathscr{C} : \sigma \neq \sigma' \text{ and } \sigma, \sigma' \text{ are tangent}\}| \lesssim_{\epsilon} B^{1/n - \epsilon} N^{2 - 1/n + \epsilon}. \]
\end{theorem}
To the best of our knowledge, the result is new in dimensions $n>3$. Observe that when $n=3$ the exponent of $N$ is $2 - 1/3 = 5/3$, which does not match the optimal exponent $3/2$; it is most unlikely that the exponent $2 - 1/n$ is asymptotically sharp in any dimension.\par
The proof combines techniques from different papers. We employ the lifting technique of \cite{EllenbergSolymosiZahl} to translate the tangency problem into a problem of incidence geometry for varieties in a space of larger dimension. To tackle the resulting incidence problem we employ the polynomial partitioning argument with bounded-degree polynomials of Solymosi and Tao (\cite{SolymosiTao}); however, instead of the standard polynomial partitioning lemma for points, we use polynomial partitioning for varieties, as developed by Guth in \cite{Guth}. The algebraic contribution to the number of tangencies is dealt with via a recursive argument.
\begin{remark}
The argument that will be given to prove Theorem \ref{main_theorem} can easily be repurposed to yield yet another proof of the $O(N^{3/2})$ bound for the number of circle tangencies in non-degenerate collections. However, since the proof relies on polynomial partitioning, it does not extend to fields other than $\mathbb{R}$.
\end{remark}
The conditional dependence of the bound in Theorem \ref{main_theorem} on the external parameter $B$ is clearly undesirable, and it is natural to look for ways to remove it. In Section 10.8 of \cite{GuthBook} an argument is presented that allows such a removal in the context of the problem of counting the number of $r$-rich intersection points for a collection of $N$ lines in $\mathbb{R}^3$. In that case the parameter $B$ controls the number of lines contained in any low-degree algebraic surface. In our context it appears that the argument of \cite{GuthBook} cannot be repeated, the main obstacle being that the low-degree varieties that contain many points of tangency are subsets of $\mathbb{R}^n$, whereas the polynomial method we use takes places in $\mathbb{R}^{2n-1}$ after the lifting procedure. While it might be possible to alter the notion of non-degeneracy we employ in order to compensate for this fact, we refrain from doing so as the resulting condition would appear somewhat unnatural.
\subsection*{Notation}
We use $\mathbb{R}[X_1,\ldots, X_m]$ to denote the ring of polynomials in $m$ variables with real coefficients; for $P$ in such a ring, we denote by $Z(P)$ the (real) zero set of $P$, that is $Z(P):= \{ x \in \mathbb{R}^m : P(x) = 0\}$. The complex zero set of $P$ is instead denoted $Z_{\mathbb{C}}(P):= \{ z \in \mathbb{C}^m : P(z) = 0\}$. Given polynomials $P_1,\ldots, P_k$ we write $Z(P_1,\ldots, P_k) = \bigcap_{j=1}^{k} Z(P_j)$.
\section{Lifting and geometric preliminaries}
The following is a straightforward generalisation of the lifting procedure in \cite{EllenbergSolymosiZahl}. We will use $X$ to denote $(X_1,\ldots, X_n)$ and $Y$ to denote $(Y_1, \ldots, Y_{n-1})$, for shortness. For $x = (x_1, \ldots, x_n) \in \mathbb{R}^n$, $r>0$ denote by $S_{x,r}$ the polynomial
\[ S_{x,r}(X) := r^2 - \sum_{j=1}^{n} (X_j - x_j)^2, \]
so that $Z(S_{x,r})$ is the sphere of radius $r$ with centre in $x$. Furthermore, for $x$ as above, we define the polynomials $Q_{x,j} \in \mathbb{R}[X_1,\ldots, X_n, Y_1, \ldots, Y_{n-1}]$ given by
\[ Q_{x,j}(X,Y):= (X_n - x_n) Y_j - (X_j - x_j), \]
where $j \in \{1,\ldots, n-1\}$. Given a sphere $\sigma$ in $\mathbb{R}^n$, let us denote its centre by $x^\sigma$ and its radius by $r_\sigma$; then we introduce the algebraic variety $\beta[\sigma]$ in $\mathbb{R}^{2n-1}$ given by
\[ \beta[\sigma] := Z(S_{x^{\sigma}, r_{\sigma}}, Q_{x^{\sigma},1}, \ldots, Q_{x^{\sigma},n-1}) \]
(where with a little abuse of notation we are regarding $S_{x^{\sigma}, r_{\sigma}}$ as a polynomial in $(X,Y)$). The variety $\beta[\sigma]$ is the union of two disjoint irreducible components, one of which is the variety given by equation $X = x^{\sigma}$ (a ``vertical'' subspace). We denote the other irreducible component\footnote{This irreducible component is cut for example by the polynomial $r_{\sigma}^2 - (X_n - x_n^{\sigma})^2 \sum_{j=1}^{n-1} Y_j^2$.} by $\beta^{\ast}[\sigma]$; this will be our ``lift'' of $\sigma$. Observe that under the projection $(x,y) \mapsto x$ this component is mapped into $\sigma$ minus the equator $X_n = x_n^\sigma$ (that is, we must have $x_n \neq x_n^{\sigma}$). The geometric meaning of the lifting is the following: if $(x,y) \in \beta^{\ast}[\sigma]$, we have that the hyperplane of $\mathbb{R}^n$ tangent to $\sigma$ at point $x$ is given by equation $(y,1) \cdot (X-x) = 0$. Indeed, we have $x_n \neq x_n^{\sigma}$, and thus $Q_{x^{\sigma},j}(x,y)=0$ implies  $y_j = (x_j - x_j^{\sigma})/(x_n - x_n^{\sigma})$; the claim then follows from elementary geometry.\par
By rotating the whole collection $\mathscr{C}$ slightly, we may assume that whenever two spheres $\sigma_1, \sigma_2$ are tangent at $x$ we have $x_n \neq x_{n}^{\sigma_1}, x_{n}^{\sigma_2}$. As a consequence, if $\sigma_1, \sigma_2$ are tangent at point $x$, we have that the intersection $\beta^{\ast}[\sigma_1] \cap \beta^{\ast}[\sigma_2]$ contains the point $(x,y)$, where $y$ is such that $(y,1)\cdot(X-x)=0$ is the equation of the common tangent hyperplane. Viceversa, if $\beta^{\ast}[\sigma_1] \cap \beta^{\ast}[\sigma_2]$ contains a point $(x,y)$ then we see (by projecting onto the first $n$ components) that $x$ is contained in $\sigma_1 \cap \sigma_2$, and moreover the tangent hyperplane to $\sigma_1$ at point $x$ coincides with the tangent hyperplane to $\sigma_2$, also at point $x$ (since it is specified by $(y,1)$); therefore $\sigma_1$ and $\sigma_2$ are actually tangent to each other at $x$. This shows in particular that two varieties $\beta^{\ast}[\sigma_1], \beta^{\ast}[\sigma_2]$ can intersect in at most a single point (thus they satisfy a pseudoline-type axiom). Since we have imposed that no three spheres can be mutually tangent to each other at a same point, the number of tangent pairs in $\mathscr{C}$ is then equal to the number of incident varieties in the collection of lifted spheres: letting $\beta^{\ast}[\mathscr{C}]:= \{ \beta^{\ast}[\sigma] : \sigma \in \mathscr{C}\}$ denote this collection,
\begin{equation}\label{eq:tangency_to_incidence}
\begin{aligned}
&|\{ (\sigma_1,\sigma_2) \in \mathscr{C} \times \mathscr{C} : \sigma_1 \neq \sigma_2 \text{ and } \sigma_1, \sigma_2 \text{ are tangent}\}| \\
& \hspace{3em}= |\{(\beta_1,\beta_2) \in \beta^{\ast}[\mathscr{C}] : \beta_1 \neq \beta_2 \text{ and } \beta_1 \cap \beta_2 \neq \varnothing \}|.
\end{aligned}
\end{equation}
\par
We will need a further important geometric fact in the next section, which is the higher dimensional analogue of a similar fact in \cite{EllenbergSolymosiZahl}. When two varieties $\beta_1,\beta_2 \in \beta^{\ast}[\mathscr{C}]$ intersect at a point, the span of their tangent spaces will contain all ``vertical'' directions of the form $(0,v) \in \mathbb{R}^n \times \mathbb{R}^{n-1}$, where $v \neq 0$ is any element of $\mathbb{R}^{n-1}$. Indeed, first of all observe that all the points of $\beta = \beta^{\ast}[\sigma]$ are smooth, since the Jacobian matrix associated to the variety has rank $n$; thus the tangent space is well-defined at every point. Now assume that the spheres $\sigma_1,\sigma_2$ are tangent to each other at $x$ and let $y$ be such that $\{(x,y)\} = \beta^{\ast}[\sigma_1] \cap \beta^{\ast}[\sigma_2]$. We write $x = (\tilde{x}, x_n) \in \mathbb{R}^{n-1} \times \mathbb{R}$ (equivalently, $\tilde{x} = (x_1, \ldots, x_{n-1})$). Since by assumption $x_n \neq x_n^{\sigma_i}$, there exist local parametrisations of the $\sigma_i$'s around $x$ of the form
\[ (\tilde{x} + t,\, x_n + \varphi_i(t)), \qquad |t| < \delta \]
for $i \in \{1,2\}$ and some $\delta > 0$, where $\varphi_i(0) = 0$; then we see that $y = - \nabla \varphi_i(0)$, and moreover we obtain local parametrisations of the $\beta^{\ast}[\sigma_i]$'s by 
\[ (\tilde{x} + t,\, x_n + \varphi_i(t),\, - \nabla \varphi_i(t)). \]
The tangent space $T_{(x,y)}\beta^{\ast}[\sigma_i]$ is then spanned by the vectors 
\[ \bm{v}^{(i)}_j := (e_j, \, \partial_j \varphi_i(0), \, - \nabla \partial_j \varphi_i(0)) \]
for $j \in \{1,\ldots, n-1\}$. We can inspect those second derivatives by hand: taking derivatives of the equation $S_{x^{\sigma_i}, r_{\sigma_i}} = 0$ that defines the sphere $\sigma_i$, a simple calculation reveals that 
\[ \partial_j \partial_k \varphi_i(0) = - \frac{\delta_{jk} (x_n - x_n^{\sigma_i})^2 + (x_j - x_j^{\sigma_i})(x_k - x_k^{\sigma_i})}{(x_n - x_n^{\sigma_i})^3}, \]
where $\delta_{jk}$ denotes the Kronecker delta; thus 
\[ -\nabla \partial_j \varphi_i(0) = \frac{1}{x_n - x_n^{\sigma_i}}\, e_j + \frac{x_j - x_j^{\sigma_i}}{(x_n - x_n^{\sigma_i})^3} \,(\tilde{x} - \tilde{x}^{\sigma_i}). \]
Notice that for some $\lambda \not\in \{0,1\}$ we have 
\[ x - x^{\sigma_2} = \lambda (x - x^{\sigma_1}), \]
since the segments from $x$ to the centres of the spheres are parallel. The subspace 
\[ \operatorname{Span}( T_{(x,y)}\beta^{\ast}[\sigma_1],  T_{(x,y)}\beta^{\ast}[\sigma_2]) \]
contains the vectors $\bm{v}^{(1)}_j - \bm{v}^{(2)}_j$, which by the observations above are equal to 
\[ (1 - \lambda^{-1}) \Big( 0, \, 0, \, \frac{1}{x_n - x_n^{\sigma_1}}\, e_j + \frac{x_j - x_j^{\sigma_1}}{(x_n - x_n^{\sigma_1})^3} \,(\tilde{x} - \tilde{x}^{\sigma_1}) \Big). \]
To show our claim it therefore suffices to prove that the vectors in the last components above are all linearly independent of each other; equivalently, letting $z := (x_n - x_n^{\sigma_1})^{-1} (\tilde{x} - \tilde{x}^{\sigma_1})$ (a column vector), the matrix whose columns are the vectors in question is given by $(x_n - x_n^{\sigma_1})^{-1}(I + z z^{\top})$, and therefore it suffices to show that the determinant 
\[ \det( I + z z^{\top}) \]
is non-zero. However, the latter is simply $1 + \|z\|^2 \neq 0$, since 
\[ \begin{pmatrix}
I & 0 \\ -z^{\top} & 1
\end{pmatrix} \begin{pmatrix}
I & z \\ 0 & 1 + \|z\|^2
\end{pmatrix} \begin{pmatrix}
I & 0 \\ z^{\top} & 1
\end{pmatrix} = \begin{pmatrix}
I + z z^{\top} & z \\ 0 & 1
\end{pmatrix}, \]
and therefore we have shown that 
\begin{equation}\label{eq:geometric_fact}
\{0\} \times \mathbb{R}^{n-1} \subset \operatorname{Span}( T_{(x,y)}\beta^{\ast}[\sigma_1],  T_{(x,y)}\beta^{\ast}[\sigma_2]), 
\end{equation}
as claimed. 
\section{Proof of Theorem \ref{main_theorem}}
Below we will take $\epsilon>0$ fixed and $\mathscr{C}$ to be a collection of $N$ spheres as per the statement of Theorem \ref{main_theorem}. The degree $D = D(\epsilon, n)$ will be considered fixed throughout, but will be chosen at the end. We will also assume it is even, for simplicity.
\subsection{Polynomial partitioning argument}
As anticipated, we will employ the inductive polynomial partitioning argument of \cite{SolymosiTao} to control the number of tangent pairs in $\mathscr{C}$. However, we will resort to a polynomial partitioning result for varieties rather than the one for points used in \cite{SolymosiTao}. The result we are going to use is due to Guth and is as follows.
\begin{theorem}[Theorem 0.3 of \cite{Guth}]\label{thm:poly_part_varieties}
Let $\Gamma$ be a collection consisting of $N$ distinct $k$-dimensional varieties in $\mathbb{R}^d$, each defined by at most $m$ polynomial equations of degree bounded by $e$. For any $D \geq 1$ there exists a polynomial $P \in \mathbb{R}[X_1, \ldots, X_d]$ of degree at most $D$ such that each connected component of $\mathbb{R}^d \setminus Z(P)$ intersects at most $C_1 N / D^{d-k}$ varieties of $\Gamma$ (where $C_1 = C_1(d,e,m)$ is a constant).
\end{theorem}
To put the above statement into context, we remark that the number of connected components of $\mathbb{R}^d \setminus Z(P)$ (where $P$ has degree at most $D$) is bounded by $C_2\, D^d$, where $C_2 = C_2(d)$ is a dimensional constant (see \cite{Milnor} and \cite{Thom}).\par
We begin by setting up the inductive argument. By \eqref{eq:tangency_to_incidence} it will suffice to work with incidences between lifted spheres: for any collection of spheres $\mathscr{D}$, lets us define the set of incidences between the lifted collections as 
\[ I(\mathscr{D}):= \{ (\beta_1, \beta_2) \in \beta^{\ast}[\mathscr{D}] \times \beta^{\ast}[\mathscr{D}] : \beta_1 \neq \beta_2 \text{ and } \beta_1 \cap \beta_2 \neq \varnothing \}. \]
For a clean inductive argument, we introduce the quantity
\[ \Theta(B,N) := \sup |I(\mathscr{C})|, \]
where the supremum is taken over all collections of spheres $\mathscr{C}$ which are $(B,D)$-non-degenerate and of cardinality $\leq N$. Our goal is to show that $\Theta(B,N) \lesssim_{\epsilon} B^{1/n - \epsilon} N^{2 - 1/n + \epsilon}$ for all $B \leq N$, which by \eqref{eq:tangency_to_incidence} proves Theorem \ref{main_theorem}.\par 
Let $P \in \mathbb{R}[X,Y]$ be the polynomial of degree at most $D/2$ produced by Theorem \ref{thm:poly_part_varieties} for the collection $\Gamma = \beta^{\ast}[\mathscr{C}]$ (so that $d=2n-1$, $k=n-1$); thus every connected component of $\mathbb{R}^{2n-1} \setminus Z(P)$ intersects at most $C_1 2^n N/ D^{n}$ varieties in $\beta^{\ast}[\mathscr{C}]$ (and the connected components are at most $C_2(D/2)^{2n-1}$ in number). We partition the set of incidences in two, according to whether the incidences occur in $Z(P)$ or not. To this end, given a collection of spheres $\mathscr{D}$ and a polynomial $R \in \mathbb{R}[X,Y]$, define $I_{Z(R)}(\mathscr{D})$ to be the set of incidences that occur in $Z(R)$, that is
\[ I_{Z(R)}(\mathscr{D}):= \{(\beta, \beta') \in I(\mathscr{D}) : \beta \cap \beta' \in Z(R)\}; \]
our partition is then written $I(\mathscr{C}) = I_1 \cup I_2$, where
\[ \begin{aligned}
I_1 &:= I(\mathscr{C}) \setminus I_{Z(P)}(\mathscr{C}), \\
I_2 &:= I_{Z(P)}(\mathscr{C}). 
\end{aligned} \]
Let us begin by bounding $|I_1|$ inductively. Since all incidences in $I_1$ occur outside of $Z(P)$, by the pigeonhole principle there must exist a connected component $O_{\ast}$ of $\mathbb{R}^{2n-1} \setminus Z(P)$ such that 
\[ |\{ (\beta, \beta') \in I(\mathscr{C}) : \beta \cap \beta' \in O_{\ast} \}| \geq \frac{|I_1|}{C_2 (D/2)^{2n-1}}. \]
Let $\mathscr{C}_{O_{\ast}}$ be the sub-collection of spheres $\sigma \in \mathscr{C}$ such that $\beta^{\ast}[\sigma]$ intersects $O_{\ast}$; this collection has cardinality at most $C_1 2^n N/ D^n$ (by construction) and is $(B,D)$-non-degenerate, so we have 
\[ |\{ (\beta, \beta') \in I(\mathscr{C}) : \beta \cap \beta' \in O_{\ast} \}| \leq |I(\mathscr{C}_{O_{\ast}})| \leq \Theta(B,C_1 2^n N/ D^n), \] 
and therefore 
\[ |I_1| \leq C_2 (D/2)^{2n-1} \Theta(B,C_1 2^n N/ D^n). \]\par
Now let us consider $I_2$. In this case we will distinguish the incidences according to whether the incident varieties themselves are contained in $Z(P)$ or not. Define then, for a given collection of spheres $\mathscr{D}$ and a polynomial $R \in \mathbb{R}[X,Y]$, the subcollection 
\[ \mathscr{D}_{Z(R)} := \{ \sigma \in \mathscr{D} : \beta^{\ast}[\sigma] \subset Z(R)\}; \]
with this notation, the desired partition $I_2 = I_3 \cup I_4$ is given by 
\[ \begin{aligned}
I_3 &:= I_{Z(P)}(\mathscr{C}) \setminus I(\mathscr{C}_{Z(P)}), \\
I_4 &:= I(\mathscr{C}_{Z(P)}) 
\end{aligned}\]
(notice that $I(\mathscr{C}_{Z(P)}) \subset I_{Z(P)}(\mathscr{C})$ indeed).\par
To estimate $|I_3|$ we will make use of the assumption of $(B,D)$-non-degeneracy. Observe that we have 
\begin{equation}
\begin{aligned}
|I_3| &= |\{(\beta,\beta') \in I(\mathscr{C}) : \beta \cap \beta' \in Z(P), \; \beta \text{ or } \beta' \in \beta^{\ast}[\mathscr{C}\setminus\mathscr{C}_{Z(P)}]\}| \\
&= 2 \sum_{\sigma \in \mathscr{C}\setminus \mathscr{C}_Z} |\{ \beta' \in \beta[\mathscr{C}] : \beta' \neq \beta^{\ast}[\sigma] \text{ and }\beta^{\ast}[\sigma] \cap \beta' \in Z(P)\}|, 
\end{aligned}\label{eq:I_3_sum}
\end{equation}
and for each $\sigma \in \mathscr{C}\setminus\mathscr{C}_{Z(P)}$ the incidences appearing in the above summand are in bijection with the corresponding tangent pairs of the form $(\sigma,\sigma')$: more precisely, we have 
\begin{align*}
&|\{ \beta' \in \beta[\mathscr{C}] : \beta' \neq \beta^{\ast}[\sigma] \text{ and } \beta^{\ast}[\sigma] \cap \beta' \in Z(P)\}| \\
& \hspace{5em} = |\{ \sigma' \in \mathscr{C} : \sigma,\sigma' \text{ are tangent and } \beta^{\ast}[\sigma] \cap \beta^{\ast}[\sigma'] \in Z(P)\}|.
\end{align*}
For $\sigma$ fixed, the tangent pairs $(\sigma,\sigma')$ are furthermore (by assumption) in bijection with the corresponding points of tangency $\sigma \cap \sigma'$; let us therefore consider the equations that such points must satisfy. Fix $\sigma \in \mathscr{C}\setminus \mathscr{C}_{Z(P)}$ and let $\sigma' \in \mathscr{C}$ be such that $\beta^{\ast}[\sigma] \cap \beta^{\ast}[\sigma'] = \{(x,y)\} \in Z(P)$; that is, $P(x,y) = 0$. We have by construction that $x_n \neq x_n^{\sigma}$, so that, since $Q_{x^{\sigma},j}(x,y)=0$, we have $y_j = (x_j - x_j^{\sigma})/(x_n - x_n^{\sigma})$. If we introduce the polynomial $\tilde{P}_{\sigma} \in \mathbb{R}[X_1,\ldots, X_n]$ given by 
\[ \tilde{P}_{\sigma}(X) := P\bigg( X_1, \ldots, X_n, \frac{X_1 - x_1^{\sigma}}{X_n - x_n^{\sigma}}, \ldots, \frac{X_{n-1} - x_{n-1}^{\sigma}}{X_n - x_n^{\sigma}}\bigg)(X_n - x_n^{\sigma})^{\deg P}, \]
then we see that $\tilde{P}_{\sigma}(x) = 0$ (since $P(x,y)=0$). In other words, the point at which $\sigma, \sigma'$ are tangent to each other lies in $Z(\tilde{P}_{\sigma})$. Observe that we have $\sigma \not\subset Z(\tilde{P}_{\sigma})$: indeed, by assumption we have $\beta^{\ast}[\sigma] \not\subset Z(P)$, so that there is a point $(x',y') \in \beta^{\ast}[\sigma]$ such that $P(x',y') \neq 0$; but since $x'_n \neq x_n^{\sigma}$ we have again $y'_j = (x'_j - x_j^{\sigma})/(x'_n - x_n^{\sigma})$ and thus $\tilde{P}_{\sigma}(x') = P(x',y')(x'_n - x_n^{\sigma})^{\deg P} \neq 0$. Since $\sigma$ is irreducible, the intersection $\sigma \cap Z(\tilde{P}_{\sigma}) = Z(S_{x^{\sigma},r_{\sigma}},\tilde{P}_{\sigma})$ is thus a genuine subvariety of $\sigma$, and it contains all the points at which $\sigma$ is tangent to one of the $\sigma'$ as above: in particular,
\begin{align*}
& |\{ \sigma' \in \mathscr{C} : \sigma,\sigma' \text{ are tangent and } \beta^{\ast}[\sigma] \cap \beta^{\ast}[\sigma'] \in Z(P)\}| \\
& \hspace{3em}= |\{ x \in \sigma \cap Z(\tilde{P}_{\sigma}) : \sigma \text{ is tangent to } \sigma' \text{ at } x, \text{ for some } \sigma' \in \mathscr{C}\}|.
\end{align*}
As $\mathscr{C}$ is $(B,D)$-non-degenerate and $\deg \tilde{P}_{\sigma} \leq D$, the last quantity is at most $B$ (notice that the degree of $Z(S_{x^{\sigma},r_{\sigma}},\tilde{P}_{\sigma})$ is indeed at most $2D$), and therefore we have by \eqref{eq:I_3_sum} that
\[ |I_3| \leq 2B |\mathscr{C} \setminus \mathscr{C}_{Z(P)}|. \]
\subsection{Algebraic contribution}
In order to complete our bound of $|I_2|$ it remains to estimate $|I_4|$, or equivalently $|I(\mathscr{C}_{Z(P)})|$ -- we will do so iteratively. If $\sigma, \sigma' \in \mathscr{C}_{Z(P)}$ are such that $\beta^{\ast}[\sigma] \cap \beta^{\ast}[\sigma'] = \{(x,y)\}$, observe that by \eqref{eq:geometric_fact} and the fact that $\beta^{\ast}[\sigma], \beta^{\ast}[\sigma'] \subset Z(P)$ we have
\[ \{0\} \times \mathbb{R}^{n-1} \subset T_{(x,y)}Z(P); \]
as a consequence, we have that
\[ (x,y) \in Z(\partial_{Y_j} P) \]
for every $j \in \{1,\ldots, n-1\}$, or in other words, for every $j$ 
\[ I(\mathscr{C}_{Z(P)}) = I_{Z(\partial_{Y_j} P)}(\mathscr{C}_{Z(P)}). \]
Assuming that for some $j_1$ the polynomial $\partial_{Y_{j_1}} P$ is not identically vanishing, the above fact means that $I_4$ is a term analogous to $I_2$, and therefore we can approach it by a similar decomposition. In order to describe the iterative procedure, define 
\[ \mathscr{C}_1 := \mathscr{C}_{Z(P)}, \qquad P_1 := \partial_{Y_{j_1}} P, \]
(so that the above is rewritten $I(\mathscr{C}_1) = I_{Z(P_1)}(\mathscr{C}_1)$) and define iteratively 
\[ \mathscr{C}_{k+1} := (\mathscr{C}_k)_{Z(P_k)} \]
(that is, $\mathscr{C}_{k+1}$ is the collection of $\sigma \in \mathscr{C}_{k}$ such that $\beta^{\ast}[\sigma] \subset Z(P_k)$) and 
\[ P_{k+1} := \partial_{Y_{j_{k+1}}}P_k, \]
where $j_{k+1} \in \{1, \ldots, n-1\}$ is such that $\partial_{Y_{j_{k+1}}}P_k$ is not identically vanishing (if one such $j_{k+1}$ exists). Notice that $\mathscr{C}_{k+1} \subset \mathscr{C}_k$, that $\deg P_{k+1} < \deg P_k$ and that by a repetition of the above argument using \eqref{eq:geometric_fact} we have for every $k$ 
\[ I(\mathscr{C}_k) = I_{Z(P_k)}(\mathscr{C}_k). \]
We can therefore write
\[ |I(\mathscr{C}_k)| = |I(\mathscr{C}_k) \setminus I(\mathscr{C}_{k+1})| + |I(\mathscr{C}_{k+1})| \]
and since $I(\mathscr{C}_k)\setminus I(\mathscr{C}_{k+1}) = I_{Z(P_k)}(\mathscr{C}_k)\setminus I((\mathscr{C}_k)_{Z(P_k)})$ we have 
\[ |I(\mathscr{C}_k) \setminus I(\mathscr{C}_{k+1})| \leq 2B|\mathscr{C}_{k+1} \setminus \mathscr{C}_{k+1}| \]
by a repetition of the argument given for $|I_3|$. The iteration stops when we reach a $k$ such that $\partial_{Y_j}P_k$ vanishes identically for all $j \in \{1,\ldots, n-1\}$ (and this must eventually happen after at most $\deg P$ steps), in which case $P_k$ is actually independent of the $Y$ variable, that is 
\[ P_k(X,Y) = R(X) \]
for some non-vanishing polynomial $R \in \mathbb{R}[X_1,\ldots, X_n]$. In this case observe that when $\sigma \in \mathscr{C}_{k+1}$ (that is, $\sigma$ is such that $\beta^{\ast}[\sigma] \subset Z(P_k)$) we have actually $\sigma \subset Z(R)$ by projecting onto $\mathbb{R}^n$. However, since $\sigma = Z(S_{x^{\sigma},r_{\sigma}})$ and $S_{x^{\sigma},r_{\sigma}}$ is irreducible, we have by Hilbert's Nullstellensatz\footnote{Indeed, $\sigma$ is Zariski-dense in $Z_{\mathbb{C}}(S_{x^{\sigma},r_{\sigma}})$ because it contains a smooth point; thus $Z_{\mathbb{C}}(S_{x^{\sigma},r_{\sigma}}) \subset Z_{\mathbb{C}}(R)$ too.} that the polynomial $S_{x^{\sigma},r_{\sigma}}$ divides $R$, and thus $|\mathscr{C}_{k+1}| \leq (\deg R)/2 \leq D$. Using the trivial bound we have then
 \[ |I(\mathscr{C}_{k+1})| \leq D^2, \]
and so collecting all the contributions we have obtained
\begin{align*}
|I_4| &\leq 2B|\mathscr{C}_1 \setminus \mathscr{C}_2| + \ldots + 2B|\mathscr{C}_{k} \setminus \mathscr{C}_{k+1}| + D^2 \\
& \leq 2B |\mathscr{C}_1| + D^2.
\end{align*}
\subsection{Conclusion}
Summing the bounds obtained for $|I_1|,|I_3|,|I_4|$, we have shown that for a $(B,D)$-non-degenerate collection $\mathscr{C}$ of cardinality $\leq N$ we have 
\[ |I(\mathscr{C})| \leq C_2 (D/2)^{2n-1} \Theta(B,C_1 2^n N/ D^n) + 2 B N + D^2; \]
taking the supremum over all such collections $\mathscr{C}$, the left-hand side is replaced by $\Theta(B,N)$. Let us write $\Theta(N)$ in place of $\Theta(B,N)$ for shortness, and let us define $C_3 := 2^n C_1$, $C_4 := 2^{1-2n}C_2$; then we have obtained the inequality
\[ \Theta(N) \leq C_4 D^{2n-1} \Theta(C_3 D^{-n} N) + 2 B N + D^2. \]
Provided $D = D(\epsilon,n)$ is chosen sufficiently large that $C_3 D^{-n} < 1$, the inequality can be used recursively: after $k$ recursions we see that we have 
\begin{align}
\Theta(N) &\leq (C_4 D^{2n-1})^k \Theta((C_3 D^{-n})^k N) \nonumber \\
& \hspace{1em} + 2BN( 1 + (C_3 C_4 D^{n-1}) + \ldots + (C_3 C_4 D^{n-1})^{k-1}) \nonumber \\
& \hspace{1em} + D^2( 1 + (C_4 D^{2n-1}) + \ldots + (C_4 D^{2n-1})^{k-1}) \nonumber \\
& \leq (C_4 D^{2n-1})^k \Theta((C_3 D^{-n})^k N) + 2BN(C_3 C_4 D^{n-1})^{k} + D^2 (C_4 D^{2n-1})^{k}.  \label{eq:recursion_end}
\end{align} 
Let $k$ be such that 
\begin{equation}\label{eq:critical_k} 
(C_3 D^{-n})^{k+1} N < B \leq (C_3 D^{-n})^k N; 
\end{equation}
then by the trivial bound $\Theta(M) \leq M^2$ we see that have
\[ (C_4 D^{2n-1})^k \Theta((C_3 D^{-n})^k N) \leq (C_4 D^{2n-1})^k C_3^{-2} D^{2n} B^2. \]
From \eqref{eq:critical_k} we further see that $D^{nk} \leq C_3^k N B^{-1}$ and that $(C_3^{1/n} D^{-1})^k \leq C_3^{-1/n} D B^{1/n} N^{-1/n}$, and therefore 
\begin{align*}
(C_4 D^{2n-1})^k C_3^{-2} D^{2n} B^2 &\leq C_3^{-2} D^{-2n} (C_3^2 C_4 D^{-1})^k N^2  \\
& \leq C_3^{-2-1/n} D^{1-2n} (C_4 C_3^{2-1/n})^k B^{1/n} N^{2 - 1/n}. 
\end{align*}
Choosing $D$ sufficiently large that $C_4 C_3^{2-1/n} \leq C_{3}^{-\epsilon} D^{n\epsilon}$ we have again from \eqref{eq:critical_k} that $(C_4 C_3^{2-1/n})^k \leq B^{-\epsilon} N^{\epsilon}$, and therefore 
\[ (C_4 D^{2n-1})^k \Theta((C_3 D^{-n})^k N) \lesssim_{\epsilon,n} B^{1/n - \epsilon} N^{2 - 1/n + \epsilon}. \]
Similar arguments show that the two remaining terms in \eqref{eq:recursion_end} are also bounded by $O_{\epsilon,n}(B^{1/n - \epsilon} N^{2 - 1/n + \epsilon})$, and therefore the proof of Theorem \ref{main_theorem} is concluded.
\subsection*{Acknowledgements} The authors would like to thank Andrei Mustata for helpful consultations on algebraic geometry, Marina Iliopoulou for introducing them to the circle tangencies problem and Cosmin Pohoata for helpful consultations on the state of the art of the subject. 
%
%
%
% Uncomment the following two lines if you want to have a bibliography
\bibliographystyle{abbrv}
\bibliography{bibliography}

\begin{thebibliography}{10}

\bibitem{Dvir}
Z.~Dvir.
\newblock On the size of {K}akeya sets in finite fields.
\newblock {\em J. Amer. Math. Soc.}, 22(4):1093--1097, 2009.

\bibitem{EllenbergSolymosiZahl}
J.~S. Ellenberg, J.~Solymosi, and J.~Zahl.
\newblock New bounds on curve tangencies and orthogonalities.
\newblock {\em Discrete Anal.}, 2016.
\newblock Paper No. 18, 22.

\bibitem{FoxPachShefferSukZahl}
J.~Fox, J.~Pach, A.~Sheffer, A.~Suk, and J.~Zahl.
\newblock A semi-algebraic version of {Z}arankiewicz's problem.
\newblock {\em J. Eur. Math. Soc. (JEMS)}, 19(6):1785--1810, 2017.

\bibitem{Guth}
L.~Guth.
\newblock Polynomial partitioning for a set of varieties.
\newblock {\em Math. Proc. Cambridge Philos. Soc.}, 159(3):459--469, 2015.

\bibitem{GuthBook}
L.~Guth.
\newblock {\em Polynomial methods in combinatorics}, volume~64 of {\em
  University Lecture Series}.
\newblock American Mathematical Society, Providence, RI, 2016.

\bibitem{Milnor}
J.~Milnor.
\newblock On the {B}etti numbers of real varieties.
\newblock {\em Proc. Amer. Math. Soc.}, 15:275--280, 1964.

\bibitem{SharirSmorodinskyValculescudeZeeuw}
M.~Sharir, S.~Smorodinsky, C.~Valculescu, and F.~de~Zeeuw.
\newblock Distinct distances between points and lines.
\newblock {\em Comput. Geom.}, 69:2--15, 2018.

\bibitem{SolymosiTao}
J.~Solymosi and T.~Tao.
\newblock An incidence theorem in higher dimensions.
\newblock {\em Discrete Comput. Geom.}, 48(2):255--280, 2012.

\bibitem{Thom}
R.~Thom.
\newblock Sur l'homologie des vari\'{e}t\'{e}s alg\'{e}briques r\'{e}elles.
\newblock In {\em Differential and {C}ombinatorial {T}opology ({A} {S}ymposium
  in {H}onor of {M}arston {M}orse)}, pages 255--265. Princeton Univ. Press,
  Princeton, N.J., 1965.

\bibitem{Wolff}
T.~Wolff.
\newblock Recent work connected with the {K}akeya problem.
\newblock In {\em Prospects in mathematics ({P}rinceton, {NJ}, 1996)}, pages
  129--162. Amer. Math. Soc., Providence, RI, 1999.

\bibitem{Zahl}
J.~Zahl.
\newblock Sphere tangencies, line incidences and {L}ie's line-sphere
  correspondence.
\newblock {\em Math. Proc. Cambridge Philos. Soc.}, 172(2):401--421, 2022.

\end{thebibliography}

\end{document}